\magnification=\magstep1

\newcount\sec \sec=0
\input Ref.macros
\input math.macros
\input labelfig.tex
\input epsf
\forwardreferencetrue
\citationgenerationtrue
\initialeqmacro

\def\cone{{\cal C}_1}
\def\cnull{{\cal C}_0}
\def\c{{\cal C}}
\def\ri{{\cal R}_i}


\title{A stationary random graph of no growth rate
}

\author{\'Ad\'am Tim\'ar}
\bigskip

\abstract{We present a random isometry-invariant subgraph of a Cayley graph such that with probability 1 its exponential growth rate does not exist.}

At the Banff workshop ``Graphs, groups and stochastics" in 2011, Vadim
Kaimanovich asked the following question:

\procl q.kaimanovich
For an arbitrary $(\Gamma, o)$ unimodular random rooted graph, does the exponential
growth rate
$$\lim_{n\to\infty} \log |B_\Gamma (o,n)|/n$$
exist almost surely? Here $B_\Gamma (o,n)$ denotes the ball of radius $n$ around $o$
in $\Gamma$.
\endprocl

We give a negative answer to the question.

\procl t.example
There is a transitive unimodular graph $G$ with a random invariant
spanning subgraph $\Gamma$ such that for any point $o\in G$ and almost every $\Gamma$
$$\liminf_{n\to\infty} \log |B_\Gamma (o,n)|/n=0\;\; and$$
$$\limsup_{n\to\infty} \log |B_\Gamma (o,n)|/n=c$$
with $c\geq 1/6$. Moreover, there exists such a $\Gamma$ that is a spanning tree.
\endprocl

Unimodular random graphs are random rooted graphs with a certain stationarity property (see the next paragraph). 
Unimodular random graphs provide a common framework to examples such as automorphism-invariant random subgraphs of transitive unimodular graphs (e.g., Cayley graphs), augmented Galton-Watson trees, Borel equivalence relations.... Many such examples arise as the local weak limit (or Benjamini-Schramm limit) of a sequence of finite graphs; a main open question by Aldous and Lyons is whether all of them arise this way, \ref b.AL/. 

Suppose that the random rooted graph $(G,o)$ has every degree bounded by $d$ almost surely, and define a new graph $(G',o)$, by adding loops to vertices of $G$ in a way that makes it $d$-regular. Choose a neighbor $x$ of $o$ in $G'$ uniformly. If the resulting distribution on the triples $(G',o,x)$ is the same as the distribution of the triples $(G',x,o)$, then we call $(G,o)$ unimodular. The bound on the degrees can be removed, and several alternative definitions exist, see \ref b.AL/.

The exponential growth rate of an infinite graph is one of the important invariants that are defined as a limit. 
For other invariants, such as the speed of the random walk, the defining limit is known to exist almost surely for every unimodular random graph, using some submodular inequality, \ref b.AL/. Hence it may be surprising that a similar argument does not work for the existence of growth.

\proofof t.example

{\bf Construction:}
Let $\Delta$ be the triangular lattice, $L=\Z_2\wr \Z$ be the
lamplighter group, where we think about each element of $L$ as a pair
$(\xi,k)$, where $\xi\in\{0,1\}^\Z$ is the status of lamps and
$k\in\Z$ is the position of the marker. Define $G=\Delta\times L\times
K_{10}$, where $K_{10}$ is the complete graph on vertex set 
$\{1,2,\ldots ,10\}$. 

Consider critical percolation on the vertices of $\Delta$. Let $\cnull$
and $\cone$ be respectively the set of closed and open components;
$\c:=\cnull\cup\cone$. There is a naturally defined oriented tree on $\c$
as vertex set: let there be an edge pointing from $x$ to $y$ ($x\to
y$), if the cluster $y$ separates $x$ from infinity and they are
adjacent (in particular, if one of them is closed then the other has
to be open). Call the set of leaves in the tree ${\cal L}_1$, and
recursively, define ${\cal L}_i$ as the set of vertices $x$ in $T$
with the property that the longest path oriented away from $x$ has
length $i$ (measured in the number of vertices in it). 

Fix sequence $c_1\ll c_2\ll\ldots$ to be defined later. Define a random
variable $X_i$ to be uniformly chosen in $\{1,\ldots , c_i\}$
($i=1,2,\ldots$). Consider the equivalence relation $\ri$ on $\Z$
where $x$ and $y$ are equivalent if
for $\alpha_i:=\sum_{j=1}^i X_j\prod_{k=1}^{j-1} c_k$ we have $[(x-\alpha_i)/c_1\ldots c_i]=[(y-\alpha_i)/c_1\ldots c_i]$ (where $[.]$ denotes the floor function).
That is, the sequence of invariant partitions defined by the $(\ri )$ is coarser and coarser, and each $\ri$ consists of
classes of $c_1\ldots c_i$ consecutive. This gives rise to a sequence
${\cal P}_i$ of
invariant coarser and coarser partitions of $L$, where each class of
${\cal P}_i$ have size $c_1\ldots c_i 2^{c_1\ldots c_i}$: let 
points $(\xi_1, k_1)$ and $(\xi_2,k_2)$ be in the same class if $k_1$
and $k_2$ are in the same class of $\ri$ and $\xi_1$ and $\xi_2$ only
differ on this class. 

\def\si{{\cal S}_i}

Let $\si$ be the set of subgraphs of $G$ induced by the finite sets of the form
$\delta\times\sigma\times K_{10}$, where $\delta\in {\cal L}_i$ and 
$\sigma\in
{\cal P}_i$. So each element of $\si$ has the form of the product of a
percolation cluster of $\Delta$ in ${\cal L}_i$, a class of the partition 
of
${\cal P}_i$ and $K_{10}$. Let $\si ^0$ be the set of those elements in
$\si$ where the $\delta$ 
above is closed, and let $\si ^1$ be the set of those where it is
open. Hence $\si=\si ^0\cup\si ^1$; call the elements of $\si$ {\it cans}. Cans in $\si ^0$ will be called {\it
type 0}, those in $\si ^1$ are {\it type 1}. 

First, we will define the edges of $\Gamma$ that go between two
distinct cans, then those that go inside
one. Suppose that $\delta\times\sigma\times K_{10}$ and
$\delta'\times\sigma'\times K_{10}$ are in $\cup\si$, and such that 
$\sigma\subset\sigma'$ and
$\delta\to\delta'$ in $T$. Then choose a random edge between 
$\delta\times\sigma\times K_{10}$ and
$\delta'\times\sigma'\times K_{10}$ and add it to $\Gamma$. 
Do it for every such pair. This way we have defined a tree on the set $\si$ of cans.

What is left is to define edges within cans. For cans of type 1, let every edge induced in the can be in $\Gamma$. (If we wanted $\Gamma$
to be a tree, we can choose a spanning tree of
the can $C$ at this point that is geodesic with
respect to some point, i.e. for some vertex $x$ in $C$ the distance of any $y$ in $C$ is the same as in this tree. The diameter of the resulting spanning tree is still logarithmic in the size, and this is the only thing we will use later.) 
For cans $\delta\times\sigma\times K_{10}$ of type 0, choose a Hamiltonian 
path in
the can (according to some fix rule which only depends on the isomorphism
type of $\delta\times\sigma\times K_{10}$). To see that such a
Hamiltonian path always exists, choose any spanning tree in
$\delta\times\sigma$, take a depth-first walk $v_1,\ldots , v_m$ that
visits every vertex of $\delta\times\sigma$ at least once and at most
degree$(\Delta\times L)=10$ times, and then replace every copy of a vertex 
$v$ in
$v_1,\ldots , v_m$ by one or more of the vertices
$(v,1),(v,2),\ldots ,(v,10)\in \delta\times\sigma\times K_{10}$ in such a 
way 
that each of these occur exactly once in the
resulting new path. This finishes the construction of $\Gamma$. It is easy to check that the resulting $\Gamma$ is ergodic.

{\bf Verification:}
We sketch here why the upper and lower growth rates are different, before giving a rigorous proof in the rest. There are radii $R$ when $B_R(o):=B_\Gamma(o,R)$ just exits a can $C$ of type 0, in the direction of the can that separates it from infinity (meaning that $B_{R}(o)$ intersects only some finite components of $\Gamma\setminus C$, but $B_{R+1}$ intersects the infinite component as well). Then with high probability some constant proportion of $B_R(o)$ is within $C$, so it has to contain some constant proportion of the Hamiltonian path in $C$, which will imply that its radius is also close to the volume (this gives the claim about the lower growth rate). On the other hand, there are radii $R$ when $B_R(o)$ just exits a can $C$ of type 1, in which case a significant proportion of the volume of $B_R(o)$ is contained in $C$ with high probability. Since the growth of the ball in $C$ is exponential, the radius increase between entering and exiting $C$ is at most the logarithm of its size. By the fast increase of the $c_i$ we can conclude that $R$ is also about logarithmic (for this, it is enough to assume that the sum of the sizes of the smaller cans that are intersected by a minimal $o$-$C$ path is of order $\log c_i$ with high probability). This gives the statement about the upper growth.

\def\i{{\rm int}}
Now we work out the above argument in details. 
For a cluster $\delta\in {\cal C}$ let $\i \delta$ be the complement of the infinite component of $\Delta\setminus \delta$ (that is, the union of all clusters that $\delta$ separates from infinity, including $\delta$). 
Denote by $\pi_\Delta$ the projection from $G$ to $\Delta$ and $\pi_{H}$ be the projection to $L\times K_{10}=:H$. From now on, condition on the event $E_0:=\{o$ is in a can $C\in {\cal S}_1\}$.
Let $C=C_1,C_2,\ldots$ be consecutive cans that an infinite simple path from $o$ visits (in particular, $\pi_\Delta (C_1), \pi_\Delta (C_2),\ldots$ is a simple path in $T$ oriented from $\pi_\Delta (C_1)$ to infinity). Note that by definition if $i<i'$ then $C_{i'}$ separates $C_i$ from infinity in $\Gamma$. We will now specify how the $c_i$'s have to be chosen. Let $c_1=1$. For $i>1$,
let $c_i$ be a number such that with probability at least $1-2^{-i}$ we have 
$$\log c_i\geq 10 c_1\ldots c_{i-1}2^{c_1\ldots c_{i-1}}|\i \pi_\Delta (C_i)|
\label e.choice
.$$
Let $E_i$ be the event that \ref e.choice/ holds. 
Let $e$ be the edge connecting $C_i$ to $C_{i+1}$ in $\Gamma$, and let $R$ be such that the ball of radius $R$ around $o$ in $\Gamma$ contains exactly one endpoint of $e$ (which implies that this ball is contained in the finite component of $\Gamma\setminus e$). The finite component of $\Gamma\setminus e$ 
arises as a union of cans $C'$, namely, it is $\cup\{C'\, :\,\pi_H(C')\subset \pi_H (C_i), \pi_\Delta (C')\subset \i (C_i)\}$. Hence, conditioned on $E_i$ we have the following upper bound for the total volume of the finite component of $\Gamma\setminus e$ 
$$|\cup\{C'\, :\,\pi_H(C')\subset \pi_H (C_i), \pi_\Delta (C')\subset \i 
(C_i)\}|=|\i\pi_\Delta (C_i)||\pi_H (C_i)|$$
$$=|\i\pi_\Delta (C_i)| 10
c_1\ldots c_i 2^{c_1\ldots c_i} \leq 10 c_i^2 2^{c_i},
\label e.2
$$ 
by \ref e.choice/.
On the other hand, 
$$|C_i|= |\pi_\Delta (C_i)| c_1\ldots c_i 2^{c_1\ldots c_i}\geq  c_1\ldots 
c_i 2^{c_1\ldots c_i}.
\label e.3
$$
Now, in the case that $C_i$ is of type 0, the restriction $C_i|_\Gamma$ is a path $P$. The edge connecting $C_i$ to $C_{i+1}$ was chosen uniformly in the construction of $\Gamma$, thus, conditioned on $P$, the endpoint of 
$e$ can be any point of $P$ with equal probability. Consequently, with probability $\geq 1/2$ the ball $B_o(R)$ contains at least $1/4$ of the points in $P$, whose length hence gives a lower bound for $R$. We conclude that conditioned on $E_0$ and $E_i$, with probability at least $1/2$ there is an $R$ such that  $B_R(o)$ has radius $R\geq c_1\ldots c_i 2^{c_1\ldots c_i}/4$ (by \ref e.3/) and has volume $\leq 10 c_i^2 2^{c_i}$ (by \ref e.2/). We obtain that
$\log |B_o(R)|/R$ gets arbitrarily close to 0 with probability at least $1/2$ (and hence, by ergodicity, with probability 1). Obviously this remains true if we do not condition on $o\in C\in {\cal S}_1$.

Suppose finally, that $C_i$ is of type 1. Then the diameter of $\Gamma|_{C_i}$ is logarithmic to its size, so conditioned on $E_i$ we have $R\leq |C_1|+\ldots +|C_{i-1}|+\log |C_i|\leq 10 c_1\ldots c_{i-1}2^{c_1\ldots c_{i-1}}|\i \pi_\Delta (C_i)|+\log 
|C_i|\leq 2\log |C_i|$
with probability at least $1/2$, using \ref e.choice/. On the other hand 
with probability at least $1/2$ the volume of $B_o(R)\cap C_i$ is at least 
$|C_i|^{1/3}$, since the choice of $e$ is uniform. We obtain
$$\log |B_o(R)|/R>1/6,$$
which holds for infinitely many $R$ almost surely (by an argument similar to the end of last paragraph).
\Qed

\medbreak
\noindent {\bf Acknowledgements.}\enspace
I thank Lewis Bowen, Russ Lyons and G\'abor Pete for valuable discussions.

Research supported by Sinergia grant CRSI22-130435 of the Swiss 
National Science foundation and by
MTA Renyi "Lendulet" Groups and Graphs Research Group.

\startbib

\bibitem[AL]{AL} Aldous, D. \and Lyons, R. (2007) Processes on unimodular 
random
networks {\it Electr. Commun. Prob.} {\bf 12}, 1454-1508.

\bibitem[12]{LP} R. Lyons and Y. Peres, {\it Probability on Trees and
Networks}, in preparation.
$\;\;\;\;\;$ http://mypage.iu.edu/$\sim$rdlyons/prbtree/prbtree.html

\endbib



\bye